\documentclass[11pt,reqno]{amsart}
\usepackage{amsmath,amssymb,amsthm,mathrsfs,dsfont}
\usepackage{CJK}
\usepackage{amsmath}
\usepackage{dsfont}
\usepackage{mathrsfs}
\usepackage{amsmath,amssymb}
\usepackage{amsfonts}
\usepackage{hyperref}
\usepackage{amsthm}
\usepackage{graphicx}
\usepackage{subfigure}
\usepackage{xcolor}
\usepackage{overpic}

\newtheorem{theorem}{Theorem}[section]
\newtheorem{lemma}{Lemma}[section]
\newtheorem{definition}{Definition}[section]
\newtheorem{proposition}{Proposition}[section]
\newtheorem{corollary}{Corollary}[section]\numberwithin{equation}{section}
\newtheorem{remark}{Remark}[section]

\newfont{\bb}{msbm10 at 12pt}

\newcommand{\bal}{\begin{aligned}}      \newcommand{\eal}{\end{aligned}}
\newcommand{\ba}{\begin{array}}      \newcommand{\ea}{\end{array}}
\newcommand{\bc}{\begin{center}}     \newcommand{\ec}{\end{center}}
\newcommand{\be}{\begin{enumerate}}  \newcommand{\ee}{\end{enumerate}}
\newcommand{\beq}{\begin{eqnarray}}  \newcommand{\eeq}{\end{eqnarray}}
\newcommand{\beQ}{\begin{eqnarray*}} \newcommand{\eeQ}{\end{eqnarray*}}
\newcommand{\bi}{\begin{itemize}}    \newcommand{\ei}{\end{itemize}}
\newcommand{\bt}{\begin{tabular}}    \newcommand{\et}{\end{tabular}}
\newcommand{\bdm}{\begin{displaymath}} \newcommand{\edm}{\end{displaymath}}

\def\qed{\hfill{Q.E.D.}\smallskip}

\begin{document}

\title{\bf Combinatorial Ricci flow on compact 3-manifolds with boundary}
\author{Xu Xu}
\maketitle

\begin{abstract}
Combinatorial Ricci flow on an ideally triangulated compact 3-manifold with boundary was introduced by Luo
as a 3-dimensional analog of Chow-Luo's combinatorial Ricci flow on a triangulated surface
and conjectured to find algorithmically the
complete hyperbolic metric on the compact 3-manifold with totally geodesic boundary.
In this paper, we prove Luo's conjecture affirmatively
by extending the combinatorial Ricci flow through the singularities of the flow
if the ideally triangulated compact 3-manifold with boundary admits such a metric.
\end{abstract}


\section{Introduction}
Motivated by Chow-Luo's combinatorial Ricci flow on triangulated surfaces \cite{CL}
and Hamilton's Ricci flow on three dimensional Riemannian manifolds \cite{H},
Luo \cite{L2} introduced combinatorial Ricci flow for hyper-ideal polyhedral metrics on ideally triangulated
compact 3-manifolds with boundary consisting of
surfaces of negative Euler characteristic.
By Moise \cite{Mo}, every compact 3-manifold can be ideally triangulated.
It is conceivable that one could use the combinatorial Ricci flow to give a new proof of Thurston's geometrization theorem for these 3-manifolds.
Luo \cite{L2} conjectured that the combinatorial Ricci flow could be used to find algorithmically the
complete hyperbolic metric with totally geodesic boundary on a compact 3-manifold with boundary.
In this paper, we prove Luo's conjecture affirmatively
if such a metric exists on the ideally triangulated compact 3-manifold with boundary.

The main tool used in the proof is the extension of dihedral angles of a hyper-ideal tetrahedron introduced by Luo-Yang \cite{LY}.
For the combinatorial Ricci flow on an ideally triangulated 3-manifold with boundary,
the hyper-ideal tetrahedra may degenerate along the flow,
which corresponds to the singularities of the flow and brings the main difficulty for the analysis of the long time behavior of the flow.
We overcome this difficulty by extending the flow through the singularities using Luo-Yang's extension.
It is shown that the solution of combinatorial Ricci flow can be uniquely extended for all time in this way.
Combining with Luo-Yang's $C^1$ smooth convex extension of co-volume function, we prove that
for any initial hyper-ideal polyhedral metric on an ideally triangulated compact 3-manifold with boundary,
the extended solution of combinatorial Ricci
flow converges exponentially fast to a complete hyperbolic metric
with totally geodesic boundary if such a metric exists on the ideally triangulated compact 3-manifold with boundary,
which confirms Luo's conjecture.

We state the results more precisely as follows.
Suppose $(M, \mathcal{T})$ is an ideally triangulated compact 3-manifold with boundary.
$E, T$ represent the sets of edges and tetrahedra in the triangulation $\mathcal{T}$ respectively.
Replacing each tetrahedron (truncated tetrahedron more precisely) in $T$ by a hyper-ideal tetrahedron
and gluing them isometrically along the hexagonal faces, we obtain a hyperbolic cone manifold with boundary.
Here a hyper-ideal tetrahedron is a compact convex polyhedron in $\mathbb{H}^3$ that is diffeomorphic to a truncated tetrahedron in $\mathbb{E}^3$,
which has four right-angled hyperbolic hexagonal faces and four hyperbolic triangular faces.
Any triangular face is required to be orthogonal to its three adjacent hexagonal faces.
Note that, by the cosine laws for hyperbolic triangles and right-angled hyperbolic hexagons,
the geometry of a hyper-ideal tetrahedron is
uniquely determined by the lengths of six edges given by pairwise intersections of four hexagonal faces.
Based on this observation, a hyper-ideal polyhedral metric on an ideally triangulated compact
3-manifold with boundary is defined as follows.

\begin{definition}[\cite{L2, LY}]
Suppose $(M, \mathcal{T})$ is an ideally triangulated compact $3$-manifold with boundary.
A hyper-ideal polyhedral metric on $(M, \mathcal{T})$ is a map $l: E\rightarrow (0, +\infty)$ such that
for any topological truncated tetrahedron in $T$,
$l_{ij}, l_{ik}, l_{ih}, l_{jk}, l_{jh}, l_{kh}$ form the lengths of six edges given by pairwise intersections
of four hexagonal faces in a hyper-ideal tetrahedron.
\end{definition}

$M$ can be taken as the hyperbolic cone manifold generated by gluing the hyper-ideal tetrahedra isometrically along
the hexagonal faces and $\partial M$ has a natural triangulation induced from the cell decomposition of $M$.
The boundary $\partial M$ is required to be consisting of surfaces with negative Euler characteristic.
For simplicity, we still call $\mathcal{T}$ a triangulation of $M$.
The set of the edges given by intersections of hexagonal faces is denoted by $E$ in the following if there is no confusing in the context.
Note that this is different from the edge set in the cell decomposition of $M$.
The space of hyper-ideal polyhedral metrics on $(M, \mathcal{T})$ is denoted by $\mathcal{L}(M, \mathcal{T})$,
which is an open subset of $\mathbb{R}^E_{>0}$.
The combinatorial Ricci curvature $K_{ij}: \mathcal{L}(M, \mathcal{T})\rightarrow (-\infty, 2\pi)$ at the edge $\{ij\}\in E$ is defined to be
$2\pi$ less the sum of dihedral angles at the edge.
A hyper-ideal polyhedral metric on the triangulated manifold $(M, \mathcal{T})$ with zero combinatorial Ricci curvature corresponds
to a complete hyperbolic metric on the manifold $M$ with totally geodesic boundary.

Luo \cite{L2} introduced the following combinatorial Ricci flow
\begin{equation}\label{combinatorial Ricci flow}
\frac{dl_{ij}}{dt}=K_{ij}
\end{equation}
for hyper-ideal polyhedral metrics on ideally triangulated compact
3-manifolds with boundary consisting of surfaces of negative Euler characteristic
and established some of the basic properties of the combinatorial Ricci flow (\ref{combinatorial Ricci flow}).
Luo further conjectured that the combinatorial Ricci flow (\ref{combinatorial Ricci flow}) could be used to
find algorithmically the complete hyperbolic metric on a compact 3-manifold with totally geodesic boundary.

The main purpose of this paper is to give an affirmative answer to Luo's conjecture under the condition
that there exists a hyper-ideal polyhedral metric with zero combinatorial Ricci curvature
on the ideally triangulated compact 3-manifold with boundary.
The main results are as follows.

\begin{theorem}\label{main theorem on uniqueness and long time existence}
Suppose $(M, \mathcal{T})$ is an ideally triangulated compact $3$-manifold with boundary.
For any initial hyper-ideal polyhedral metric in $\mathcal{L}(M, \mathcal{T})$,
the solution of combinatorial Ricci flow $(\ref{combinatorial Ricci flow})$ can be extended to be
a solution existing for all time.
\end{theorem}

\begin{theorem}\label{main theorem on convergence}
Suppose $(M, \mathcal{T})$ is an ideally triangulated compact $3$-manifold with boundary.
If there exists a hyper-ideal polyhedral metric $l^*\in \mathcal{L}(M, \mathcal{T})$
with zero combinatorial Ricci curvature,
the extended solution of combinatorial Ricci flow converges exponentially fast to $l^*$
for any initial hyper-ideal polyhedral metric in $\mathcal{L}(M, \mathcal{T})$.
\end{theorem}

The local convergence of combinatorial Ricci flow (\ref{combinatorial Ricci flow}) was obtained by Luo \cite{L2}.
Theorem \ref{main theorem on uniqueness and long time existence} and Theorem \ref{main theorem on convergence}
give the global convergence of the extended combinatorial Ricci flow.
Theorem \ref{main theorem on uniqueness and long time existence} and Theorem \ref{main theorem on convergence}
imply an algorithm to find the hyper-ideal
polyhedral metric with zero combinatorial Ricci curvature
if such a metric exists on the ideally triangulated compact 3-manifold with boundary $(M, \mathcal{T})$.
Note that a hyper-ideal
polyhedral metric with zero combinatorial Ricci curvature on $(M, \mathcal{T})$ corresponds to a complete
hyperbolic metric on $M$ with totally geodesic boundary,
the existence of which is necessary for an algorithm in Luo's conjecture for the triangulated manifold $(M, \mathcal{T})$.
Therefore, Theorem \ref{main theorem on uniqueness and long time existence} and Theorem \ref{main theorem on convergence}
confirm Luo's conjecture for combinatorial Ricci flow on ideally triangulated compact 3-manifolds with boundary.

Suppose $(M, \mathcal{T})$ is an ideally triangulated compact 3-manifold with boundary,
we use $d_e$ to denote the number of tetrahedra adjacent to an edge $e\in E$.
As an application of Theorem \ref{main theorem on uniqueness and long time existence} and Theorem \ref{main theorem on convergence},
we have the following result.
\begin{corollary}\label{main corollary for regular triangulation}
Suppose $(M, \mathcal{T})$ is an ideally triangulated compact $3$-manifold with boundary.
\begin{description}
  \item[(1)] If $d_e\leq 6$ for any $e\in E$,
  there exists no hyper-ideal polyhedral metric on $(M, \mathcal{T})$ with zero combinatorial Ricci curvature.
  \item[(2)] If $d_e=d_{e'}=N>6$ for any $e, e'\in E$, then
  for any initial hyper-ideal polyhedral metric in $\mathcal{L}(M, \mathcal{T})$,
  the extended solution of combinatorial Ricci flow converges exponentially fast to
  a hyper-ideal polyhedral metric in $\mathcal{L}(M, \mathcal{T})$
  with zero combinatorial Ricci curvature.
\end{description}
\end{corollary}

One can also modify the combinatorial Ricci flow (\ref{combinatorial Ricci flow})
to find hyper-ideal polyhedral metrics with prescribed combinatorial
Ricci curvatures and then use Luo-Yang's extension to study the long time behavior of the modified combinatorial Ricci flow.

There have been many work on combinatorial curvature flows on three dimensional manifolds.
See the work of Dai-Ge \cite{DG}, Ge-Hua \cite{GH}, Ge-Jiang-Shen \cite{GJS}, Ge-Ma \cite{GM}, Ge-Xu \cite{GX1,GX2,GX4},
Ge-Xu-Zhang \cite{GXZ},  Glickenstein \cite{G1,G2}  and others.
The difference between their work and ours is that we consider compact hyperbolic cone 3-manifolds with boundary
generated by isometrically gluing hyper-ideal tetrahedra.
The combinatorial Ricci flow for compact 3-manifolds with boundary
have applications in engineering fields such as shape classification, see \cite{YJLG} for example.
A work closely related to this work is the combinatorial Ricci flow for decorated hyperbolic polyhedral metrics
on cusped $3$-manifolds introduced by Yang \cite{Y} following Luo's combinatorial Ricci flow for hyper-ideal polyhedral metrics
on compact $3$-manifolds with boundary. The author \cite{X3} 
recently extended Yang's combinatorial Ricci flow on ideally triangulated cusped $3$-manifolds 
using Luo-Yang's extension \cite{LY} and proved that
the existence of a complete hyperbolic metric on a cusped 3-manifold
is equivalent to the convergence of the extended combinatorial Ricci flow, which gives a new
characterization of existence of a complete hyperbolic metric on a cusped 3-manifold  dual to Casson and Rivin's programm.
The extended combinatorial Ricci flow also provides an elegant
and effective algorithm for finding complete hyperbolic metrics on cusped 3-manifolds.
The author \cite{X4} further introduced the combinatorial Calabi flow on cusped $3$-manifolds to find 
complete hyperbolic metrics on such manifolds. The basic properties of the combinatorial Calabi flow were established 
in \cite{X4}.
There have been many work using the extension method to study the long time behavior of combinatorial curvature
flows on low dimensional manifolds. See Ge-Hua \cite{GH}, Ge-Jiang \cite{GJ0,GJ1,GJ2,GJ3}, Ge-Jiang-Shen \cite{GJS},
Ge-Xu \cite{GX3}, Gu-Guo-Luo-Sun-Wu \cite{GGLSW}, Gu-Luo-Sun-Wu \cite{GLSW}, Xu \cite{X1,X2,X3},
Zhu-Xu \cite{ZX} and others.

The paper is organized as follows.
In Section \ref{section 2}, we recall some basic properties of hyper-ideal tetrahedra;
In Section \ref{section 3}, we recall the extension of dihedral angles introduced by Luo-Yang \cite{LY};
In Section \ref{section 4}, we introduce the extension of the combinatorial Ricci flow (\ref{combinatorial Ricci flow}) and
prove a generalization of Theorem \ref{main theorem on uniqueness and long time existence};
In Section \ref{section 5}, we prove generalizations of Theorem \ref{main theorem on convergence}
and Corollary \ref{main corollary for regular triangulation}.
In Section \ref{section 6}, we study the modified combinatorial Ricci flow to find hyper-ideal
polyhedral metrics with prescribed combinatorial Ricci curvatures.
In Section \ref{section 7}, we give some remarks and propose some questions.

\section{Preliminaries on hyper-ideal tetrahedra}\label{section 2}
In this section, we recall some basic properties of
hyper-ideal tetrahedra and set up the notations used in the following of the paper.
For more details on hyper-ideal tetrahedra, please refer to \cite{B,BB,FP,F,L2,LY,S}.

\begin{figure}[!htb]
\centering
  \includegraphics[height=0.48\textwidth,width=0.6\textwidth]{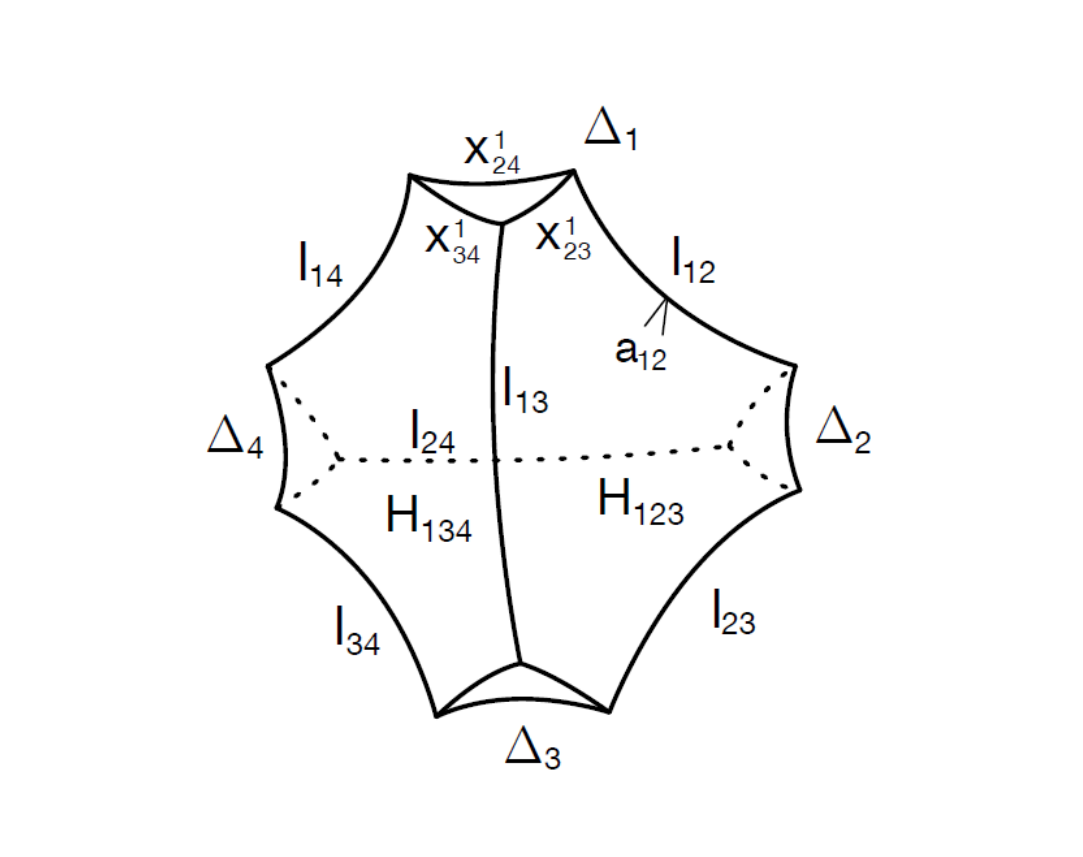}
  \caption{hyper-ideal tetrahedron (this figure is produced by Luo and Yang )}
  \label{hyper-ideal-tetrahedron}
\end{figure}

Suppose $\sigma$ is a hyper-ideal tetrahedron,
which is a compact convex polyhedron in hyperbolic space $\mathbb{H}^3$
diffeomorphic to a truncated tetrahedron in the Euclidean space $\mathbb{E}^3$
with four hexagonal faces and four triangular faces (see Figure \ref{hyper-ideal-tetrahedron}).
Denote the four triangular faces of $\sigma$ as $\Delta_i, i=1,2,3,4,$, which are hyperbolic triangles and called vertex triangles in the following.
The edge joining $\Delta_i$ and $\Delta_j$ is denoted by $e_{ij}$, the length of which is $l_{ij}$.
The hexagonal face adjacent to $e_{ij}, e_{jk}, e_{ik}$ is denoted by $H_{ijk}$, which is a right-angled hyperbolic hexagon.
The vertex triangle is required to be orthogonal to the three adjacent hexagonal faces.
The intersection of two hexagonal faces is called an edge and
the intersection of a hexagonal face and a vertex triangle is called a vertex edge.
The dihedral angle at $e_{ij}$ is the angle between the two hexagonal faces $H_{ijk}$ and $H_{ijh}$, which
is denoted by $a_{ij}$. The length of the vertex edge $\Delta_i\cap H_{ijk}$ is denoted by $x^i_{jk}$.
Let $\mathcal{L}$ be the set of vectors $(l_{12}, \cdots, l_{34})\in \mathbb{R}^6_{>0}$
such that there exists a nondegenerate hyper-ideal tetrahedron with
$l_{ij}$ as the length of the edge $\{ij\}$.

\begin{proposition}[\cite{B,BB,FP,F,L2}]\label{simply connectness of L}
$\mathcal{L}$ is a simply connected open subset of $\mathbb{R}^6_{>0}$.
\end{proposition}

For any hyper-ideal polyhedral metric $l\in \mathcal{L}$ on $\sigma$, one can define the volume function $V$ and the co-volume function
\begin{equation}\label{co-volume function for hyper-ideal tetrahedron}
F_\sigma=2V+\sum_{i<j}a_{ij}l_{ij}.
\end{equation}
By the Schl\"{a}fli formula \cite{B1,L2a,M,V}, we have
$$dF_\sigma=\sum_{i<j}a_{ij}dl_{ij},$$
which implies the matrix $\Lambda_\sigma:=[\frac{\partial a_{ij}}{\partial l_{kh}}]_{6\times 6}=Hess_l F_\sigma$ is symmetric and the 1-form $\sum_{i<j}a_{ij}dl_{ij}$
is closed on $\mathcal{L}$.
Combining with Proposition \ref{simply connectness of L}, $\int \sum_{i<j}a_{ij}dl_{ij}$ is well-defined on $\mathcal{L}$
and equals to the co-volume function $F_\sigma$ up to a constant.
Furthermore, we have the following property of the matrix $\Lambda_\sigma$.

\begin{theorem}[\cite{L2, S}]\label{positive definiteness of Lambda sigma}
For a hyper-ideal tetrahedron $\sigma$ with dihedral angle $a_{ij}$ and length $l_{ij}$
at the edge $\{ij\}$, the matrix $\Lambda_\sigma=[\frac{\partial a_{ij}}{\partial l_{kh}}]$
is symmetric and strictly positive definite on $\mathcal{L}$.
In particular, the co-volume function $F_\sigma= 2V +\sum_{i<j}a_{ij}l_{ij}$ is a locally strictly convex smooth function of
the length variables $(l_{12}, \cdots, l_{34})\in \mathcal{L}$ with $\frac{\partial F_\sigma}{\partial l_{ij}}=a_{ij}$.
\end{theorem}
Luo \cite{L2} further introduced the following function
\begin{equation}\label{energy function}
H(l)=\sum_{\sigma\in T}F_{\sigma}-2\pi\sum_{\{ij\}\in E}l_{ij},
\end{equation}
for $l\in \mathcal{L}(M, \mathcal{T})$ and proved the following result.
\begin{theorem}[\cite{L2}]\label{negative definiteness of Lambda}
The combinatorial Ricci flow (\ref{combinatorial Ricci flow}) is the negative gradient flow of the locally strictly convex
function $H$ defined on $\mathcal{L}(M, \mathcal{T})$.
Specially, the matrix $\Lambda=(\frac{\partial K_{ij}}{\partial l_{kh}})$ is
symmetric and strictly negative definite on $\mathcal{L}(M, \mathcal{T})$.
\end{theorem}

Based on Theorem \ref{positive definiteness of Lambda sigma} and Theorem \ref{negative definiteness of Lambda},
Luo \cite{L2} proved the local rigidity of hyper-ideal polyhedral metrics in $\mathcal{L}(M, \mathcal{T})$ and
the local convergence of combinatorial Ricci flow (\ref{combinatorial Ricci flow}) on ideally triangulated
compact 3-manifolds with boundary.
However, the open set $\mathcal{L}(M, \mathcal{T})$ may not be convex,
which causes the main difficulty for the global rigidity of hyper-ideal polyhedral metrics and the global
convergence of combinatorial Ricci flow (\ref{combinatorial Ricci flow}).
Luo-Yang \cite{LY} proved the global rigidity of hyper-ideal polyhedral metrics.

\section{Luo-Yang's extension of dihedral angles}\label{section 3}
To prove the global rigidity of hyper-ideal polyhedral metrics on compact 3-manifolds with boundary,
Luo-Yang \cite{LY} carefully analysed the set of degenerate hyper-ideal polyhedral metrics on a hyper-ideal tetrahedron.
Based on the characterization of the set of degenerate hyper-ideal polyhedral metrics on a hyper-ideal tetrahedron,
they extended the dihedral angles
of a hyper-ideal tetrahedron to be defined on $\mathbb{R}^6$, from which they proved the global rigidity of hyper-ideal polyhedral metrics.
We recall Luo-Yang's extension in this section, which is the main tool we use to extend the combinatorial Ricci flow (\ref{combinatorial Ricci flow}).
For more details of the extension, please refer to \cite{LY}.

For $(l_{12}, \cdots, l_{34})\in \mathbb{R}^6_{>0}$ and $\{i,j,k,h\}=\{1,2,3,4\}$, let $l_{ji}=l_{ij}$ for $i\neq j$ and let
\begin{equation*}
x^i_{jk}=\cosh^{-1}\left(\frac{\cosh l_{ij}\cosh l_{ik}+\cosh l_{jk}}{\sinh l_{ij}\sinh l_{ik}}\right)
\end{equation*}
and
\begin{equation*}
\phi^i_{kh}=\frac{\cosh x^i_{jk}\cosh x^i_{jh}-\cosh x^i_{kh}}{\sinh x^i_{jk}\sinh x^i_{jh}}.
\end{equation*}
It is proved that $\phi^i_{kh}=\phi^j_{kh}$ for any $l\in \mathbb{R}^6_{>0}$ (\cite{LY}, Lemma 4.3).
Therefore, one can define $\phi_{ij}: \mathbb{R}^6_{>0}\rightarrow \mathbb{R}$ by $\phi_{ij}(l)=\phi^i_{kh}(l)$.
If $(l_{12}, \cdots, l_{34})\in \mathcal{L}\subset \mathbb{R}^6_{>0}$ are lengths of edges in a hyper-ideal tetrahedron,
then $\phi_{ij}=\cos a_{ij}$, where $a_{ij}$ is the dihedral angle along the edge $\{ij\}$.

For $\{i,j\}\subset \{1,2,3,4\}$, let $\Omega^{\pm}_{ij}=\{l\in \mathbb{R}^6_{>0}|\pm\phi_{ij}(l)\geq 1\}$
and $X^{\pm}_{ij}=\{l\in \mathbb{R}^6_{>0}|\phi_{ij}(l)= \pm 1\}$.
Luo-Yang proved the following property of $\Omega^{\pm}_{ij}$ and $X^{\pm}_{ij}$.
\begin{lemma}[\cite{LY}, Lemma 4.6]\label{degenerate set relationship}
For $\{i,j,k,h\}=\{1,2,3,4\}$, we have
\begin{description}
  \item[(1)] $\Omega^{-}_{ij}\cap \Omega^{-}_{ik}=\emptyset$;
  \item[(2)] $\Omega^{-}_{ij}=\Omega^{-}_{kh}$ and $\Omega^{+}_{ij}=\Omega^{-}_{ik}\cup \Omega^{-}_{ih}$;
  \item[(3)] $X^{-}_{ij}=X^{-}_{kh}$ and $X^{+}_{ij}=X^{-}_{ik}\cup X^{-}_{ih}$.
\end{description}
\end{lemma}
Set $\Omega_1=\Omega^{-}_{12}$, $\Omega_2=\Omega^{-}_{13}$, $\Omega_3=\Omega^{-}_{14}$,
$X_1=X^{-}_{12}$, $X_2=X^{-}_{13}$, $X_3=X^{-}_{14}$.
As a corollary of Lemma \ref{degenerate set relationship},
Luo-Yang obtained the following characterization of $\mathcal{L}$ and $\mathbb{R}^6_{>0}\setminus \mathcal{L}$.

\begin{proposition}[\cite{LY}, Proposition 4.5]
Let $\partial \mathcal{L}$ be the frontier of $\mathcal{L}$ in $\mathbb{R}^6_{>0}$. Then $\partial \mathcal{L}=X_1\sqcup X_2 \sqcup X_3$,
where $X_i, i=1,2,3$ is a real analytic codimension-1 submanifold of $\mathbb{R}^6_{>0}$.
The complement $\mathbb{R}^6_{>0}\setminus \mathcal{L}$ is a disjoint union of three manifolds $\Omega_i$ with
boundary so that $\Omega_i\cap \partial \mathcal{L}=X_i, i=1,2,3$.
\end{proposition}

Based on this characterization, Luo-Yang \cite{LY} introduced the following extension of the dihedral angle $a_{ij}$.
For $i\neq j$, define $a_{ij}|_{\overline{\Omega^+_{ij}}}=0$ and $a_{ij}|_{\overline{\Omega^-_{ij}}}=\pi$,
where $\overline{\Omega^{\pm}_{ij}}$ is the closure of $\Omega^{\pm}_{ij}$ in $\mathbb{R}^6_{\geq 0}$.
Then the dihedral angle $a_{ij}: \mathcal{L}\rightarrow \mathbb{R}$ is continuously extended to be defined on $\mathbb{R}^6_{\geq 0}$.

Luo-Yang \cite{LY} further extended $a_{ij}$ to be defined on $\mathbb{R}^6$ as follows.
For each $l=(l_{12}, \cdots, l_{34})\in \mathbb{R}^6$, let $l^+=(l^+_{12}, \cdots, l^+_{34})\in \mathbb{R}^6_{\geq 0}$, where
$l_{ij}^+=\max\{0, l_{ij}\}$. For any $l\in \mathbb{R}^6$, set
$$\widetilde{a}_{ij}(l)=a_{ij}(l^+),$$
then $\widetilde{a}_{ij}$ is a continuous extension of $a_{ij}$ defined on $\mathbb{R}^6$.

Based on the extension theory of closed 1-forms and convex functions established in \cite{BPS,L3},
Luo-Yang \cite{LY} obtained the following extension of
the co-volume function of $F_\sigma$.
\begin{theorem}[\cite{LY}, Corollary 4.12]\label{extended co-volume}
The function $\widetilde{F}_\sigma$ defined by
\begin{equation}
\widetilde{F}_\sigma=\int_{(0,\cdots, 0)}^{l}\sum_{i\neq j}\widetilde{a}_{ij}(l)dl_{ij}+F_{\sigma}(0,\cdots, 0)
\end{equation}
is a well-defined $C^1$-smooth convex function defined on $\mathbb{R}^6$
extending the co-volume function $F_\sigma=2V+\sum_{i<j} a_{ij}l_{ij}$
defined on $\mathcal{L}\subset \mathbb{R}^6_{>0}$,
where $F_{\sigma}(0,\cdots, 0)=16 \Psi(\pi/4)$ with $\Psi$ being the Lobachevsky function.
\end{theorem}

\begin{remark}\label{remark non C1 of a_ij}
By direct calculations, Luo-Yang \cite{LY} proved that $\frac{\partial \phi_{ij}}{\partial l_{kl}}\neq 0$ for $l\in X_i$.
This implies $\widetilde{a}_{ij}$ is $C^0$ smooth and not $C^1$ smooth.
As a corollary, $\widetilde{F}_\sigma$ is not $C^2$ smooth on $\mathbb{R}^6$.
\end{remark}

Using the extension of the co-volume function $F_\sigma$ in Theorem \ref{extended co-volume},
Luo-Yang \cite{LY} proved the following global rigidity of hyper-ideal polyhedral metrics.

\begin{theorem}[\cite{LY}, Theorem 1.2 (b)]\label{LY global rigidity}
For any ideally triangulated compact $3$-manifold with boundary $(M, \mathcal{T})$, a hyper-ideal polyhedral metric on $(M, \mathcal{T})$
is determined by its combinatorial Ricci curvature, that is,
the curvature map $K: \mathcal{L}(M, \mathcal{T})\rightarrow \mathbb{R}^E$ is injective.
\end{theorem}

\begin{remark}
Luo-Yang \cite{LY} proved the global rigidity for hyper-ideal polyhedral metrics on triangulated compact pseudo $3$-manifolds,
which are generalizations of triangulated compact $3$-manifolds.
Please refer to \cite{LY} for more details.
\end{remark}

Using the extension $\widetilde{a}_{ij}$ of the dihedral angle $a_{ij}$,
Luo-Yang \cite{LY} extended the combinatorial Ricci curvature to be defined
on $\mathbb{R}^E$ by
\begin{equation}\label{generalized curvature}
\widetilde{K}_{ij}=2\pi-\sum_{\sigma\in T}\widetilde{a}_{ij},
\end{equation}
where the summation is taken with respect to the tetrahedra adjacent to the edge $\{ij\}\in E$.
In the following, $l\in \mathbb{R}^E$ is called as a generalized hyper-ideal polyhedral metric
and $\widetilde{K}_{ij}(l)$ is called the generalized combinatorial Ricci curvature
at the edge $\{ij\}\in E$ for the generalized hyper-ideal polyhedral metric $l\in \mathbb{R}^E$.

For the proof of the main results in this paper, we give the following result on rigidity of
the generalized hyper-ideal polyhedral metrics, which is a slight
generalization of Luo-Yang's global rigidity in Theorem \ref{LY global rigidity}.

\begin{theorem}\label{generalization of Luo-Yang rigidity}
For an ideally triangulated compact $3$-manifold $(M, \mathcal{T})$,
suppose $\overline{K}=K(\overline{l})$ for some $\overline{l}\in \mathcal{L}(M, \mathcal{T})$.
If there exists a generalized hyper-ideal polyhedral metric $l^*\in \mathbb{R}^E$ such that $\widetilde{K}(l^*)=\overline{K}$, then $l^*=\overline{l}$.
\end{theorem}
\proof
The proof follows essentially Luo-Yang's proof of Theorem \ref{LY global rigidity} in \cite{LY}.
For completeness, we give the proof here.
Suppose $l^*\neq \overline{l}$. Define
\begin{equation}
\widetilde{F}(l)=\sum_{\sigma\in T} \widetilde{F}_\sigma(l)
\end{equation}
for $l\in \mathbb{R}^E$. Then $\widetilde{F}$ is $C^1$ smooth and convex on $\mathbb{R}^E$
and smooth and locally strictly convex on $\mathcal{L}(M, \mathcal{T})$ with
$\nabla_{l_{ij}}\widetilde{F}=2\pi-\widetilde{K}_{ij}$ by Theorem \ref{positive definiteness of Lambda sigma} and Theorem \ref{extended co-volume}.
Define
$$f(t)=\widetilde{F}(tl^*+(1-t)\overline{l})$$
for $t\in [0,1]$.
Then $f(t)$ is a $C^1$ smooth convex function on $[0,1]$ and a smooth strictly convex function in a neighborhood $[0, \epsilon)$
for some positive number $\epsilon<1$,
which implies $f'(t)$ is a monotone nondecreasing function on $[0,1]$ and strictly nondecreasing on $[0, \epsilon)$.

Note that
\begin{equation*}
f'(t)=\sum_{\{ij\}\in E}\left[2\pi-\widetilde{K}_{ij}(tl^*+(1-t)\overline{l})\right]\cdot (l^*_{ij}-\overline{l}_{ij}),
\end{equation*}
we have $f'(0)=f'(1)$ by $\widetilde{K}(l^*)=K(\overline{l})$.
This implies $f'(t)$ is a constant on $[0, 1]$, which contradicts that $f'(t)$ is strictly nondecreasing on $[0, \epsilon)$.
\qed

\section{Uniqueness and long time existence of the extended combinatorial Ricci flow}\label{section 4}
Using the extended combinatorial Ricci curvature $\widetilde{K}_{ij}$ defined by (\ref{generalized curvature}), we can define
a new flow as follows for $l\in \mathbb{R}^E$
\begin{equation}\label{extended flow}
\frac{dl_{ij}}{dt}=\widetilde{K}_{ij}.
\end{equation}
Note that $\widetilde{K}_{ij}$ is continuous, therefore the solution of the new flow (\ref{extended flow}) exists for short time
by the ODE theory (See \cite{Hart} Theorem 2.1 for example).
However, $\widetilde{K}_{ij}$ is not locally Lipschitz continuous by Remark \ref{remark non C1 of a_ij},
the uniqueness of the solution for the new flow (\ref{extended flow}) can not be derived directly from the standard ODE theory.
For the flow (\ref{extended flow}), we still have the uniqueness of the solution.

\begin{theorem}\label{contex uniqueness}
The solution of the flow (\ref{extended flow}) is unique for any initial generalized hyper-ideal polyhedral metric $l_0\in \mathbb{R}^E$.
\end{theorem}
\proof
The idea of the proof comes from Ge-Hua \cite{GH}, but the proof is simpler for our case.
We present the proof here for completeness.
Suppose $l_1(t)$ and $l_{2}(t)$ are two solutions of the flow (\ref{extended flow}) on $[0, T]$
for the same initial value $l_0\in \mathbb{R}^E$, where $T$ is some positive constant.
To prove the uniqueness of the solution, we just need to prove $l_1(t)=l_2(t)$ for any $t\in [0, T]$.

By Theorem \ref{extended co-volume},
the energy function $H$ defined by (\ref{energy function}) could be extended to be
a $C^1$ smooth convex function
\begin{equation}\label{extension of H}
\widetilde{H}(l)=\sum_{\sigma\in T}\widetilde{F}_{\sigma}-2\pi\sum_{\{ij\}\in E}l_{ij}
\end{equation}
defined on $\mathbb{R}^E$ with $\nabla \widetilde{H}=-\widetilde{K}$.
We claim that
$$\left(\nabla \widetilde{H}(l_1)-\nabla \widetilde{H}(l_2)\right)\cdot (l_1-l_2)\geq 0,$$
which is equivalent to
\begin{equation}\label{convex inequality}
\left(\widetilde{K}(l_1)-\widetilde{K}(l_2)\right)\cdot (l_1-l_2)\leq 0,
\end{equation}
for any $l_1, l_2\in \mathbb{R}^E$.

To prove the claim, take $\varphi_\epsilon(x)=\frac{1}{\epsilon^{|E|}}\varphi(\frac{x}{\epsilon})$, $x\in \mathbb{R}^{E}$,
as the standard mollifier with
\begin{equation*}
\begin{aligned}
\varphi(x)=\left\{
             \begin{array}{ll}
               Ce^{\frac{1}{1-|x|^2}}, & \hbox{$|x|<1$,} \\
               0, & \hbox{$|x|\geq 1$,}
             \end{array}
           \right.
\end{aligned}
\end{equation*}
where $C$ is chosen to satisfy $\int_{\mathbb{R}^E}\varphi(x)dx=1$.
The $\epsilon$-mollifier $\widetilde{H}_\epsilon$ of $\widetilde{H}$ is defined to be
$$\widetilde{H}_\epsilon(l)=\widetilde{H}*\varphi_\epsilon(l)$$ for $l\in \mathbb{R}^E$.
By the convexity and $C^1$ smoothness of $\widetilde{H}$, $\widetilde{H}_\epsilon$ is a $C^\infty$ smooth convex function on $\mathbb{R}^E$
and $\widetilde{H}_\epsilon\rightarrow \widetilde{H}$ in $C^1_{loc}(\mathbb{R}^E)$.
By the $C^\infty$ smoothness and convexity of $\widetilde{H}_\epsilon$, we have
$$\left(\nabla \widetilde{H}_\epsilon(l_1)-\nabla \widetilde{H}_\epsilon(l_2)\right)\cdot (l_1-l_2)\geq 0$$
for any $l_1, l_2\in \mathbb{R}^E$.
Let $\epsilon\rightarrow 0$,  we have
$$\left(\nabla \widetilde{H}(l_1)-\nabla \widetilde{H}(l_2)\right)\cdot (l_1-l_2)\geq 0,$$
which completes the proof of the claim.

Set $f(t)=||l_1(t)-l_2(t)||^2$.
Then we have $f(t)\geq 0$ for $t\in [0,T]$ with $f(0)=0$ by the initial condition $l_1(0)=l_2(0)$.
Furthermore,
$$\frac{df}{dt}=(\frac{dl_1}{dt}-\frac{dl_2}{dt})\cdot (l_1(t)-l_2(t))
=(\widetilde{K}(l_1(t))-\widetilde{K}(l_2(t)))\cdot (l_1(t)-l_2(t))\leq 0$$
by (\ref{convex inequality}),
which implies $f(t)\equiv 0$ for $t\in [0,T]$.
Therefore, $l_1(t)=l_2(t)$ for any $t\in [0,T]$, from which the uniqueness follows.
\qed

\begin{remark}
Although the solution of the extended flow (\ref{extended flow})
is unique, there may exist some other different extensions of the solution of 
combinatorial Ricci flow (\ref{combinatorial Ricci flow}) on $\mathbb{R}^E$.
The key point is that the solution of the extended flow (\ref{extended flow}) depends on the extension of the combinatorial Ricci 
curvature.
For example, one can also extend the combinatorial Ricci curvature by symmetry to $\mathbb{R}^E$, 
which is different from the extension used here.
The author thanks Tian Yang for pointing this out to the author.
\end{remark}

As a corollary of Theorem \ref{contex uniqueness},
we have the following result which shows that the solution of the new flow (\ref{extended flow}) extends the solution of the
combinatorial Ricci flow (\ref{combinatorial Ricci flow}) for any initial hyper-ideal polyhedral metric in $\mathcal{L}(M, \mathcal{T})$.

\begin{corollary}\label{context unique extension}
For any initial hyper-ideal polyhedral metric $l(0)\in \mathcal{L}(M, \mathcal{T})$,
denote the solutions of the combinatorial Ricci flow (\ref{combinatorial Ricci flow})
and the flow (\ref{extended flow})  as $l(t)$ and $\widetilde{l}(t)$ respectively.
Then
$\widetilde{l}(t)=l(t)$
whenever $l(t)$ exists.
\end{corollary}

We call the flow (\ref{extended flow}) as the extended combinatorial Ricci flow in the following.

\begin{theorem}\label{context long time existence}
The solution of extended combinatorial Ricci flow (\ref{extended flow}) exists for all time
for any initial generalized hyper-ideal polyhedral metric $l(0)\in \mathbb{R}^E$.
\end{theorem}
\proof
By definition, the extension $\widetilde{a}_{ij}$ of the dihedral angle function $a_{ij}$ is bounded by $\pi$.
By the finiteness of the triangulation $\mathcal{T}$,
the extended combinatorial Ricci curvature $\widetilde{K}$ is uniformly bounded by some constant $C$ depending on
the triangulation $\mathcal{T}$ of $M$.
If $l(t)$ is the solution of extended combinatorial Ricci flow (\ref{extended flow}), then we have
$$\left|\frac{dl_{ij}}{dt}\right|\leq C$$
for any $\{ij\}\in E$,
which implies $|l_{ij}(t)|\leq |l_{ij}(0)|+Ct$
for any edge $\{ij\}\in E$ and $t\in [0, +\infty)$.
Therefore, the solution of extended combinatorial Ricci flow (\ref{extended flow}) does not go to infinity in finite time,
which implies the solution of extended combinatorial Ricci flow (\ref{extended flow})
exists for all time by ODE theory.
\qed

Theorem \ref{contex uniqueness}, Corollary \ref{context unique extension} and Theorem \ref{context long time existence}
together imply Theorem \ref{main theorem on uniqueness and long time existence}.

\section{Convergence of the extended combinatorial Ricci flow}\label{section 5}
\begin{theorem}
Suppose the solution $l(t)$ of extended combinatorial Ricci flow $(\ref{extended flow})$
converges to a generalized hyper-ideal polyhedral metric $l^*\in \mathbb{R}^E$, then
$\widetilde{K}_{ij}(l^*)=0$ for any edge $\{ij\}\in E$.
\end{theorem}
\proof
By the continuity of the generalized combinatorial Ricci curvature $\widetilde{K}_{ij}(l)$ in $l\in \mathbb{R}^E$,
we have
\begin{equation}\label{equ limit K_ij}
\widetilde{K}_{ij}(l^*)=\lim_{t\rightarrow +\infty}\widetilde{K}_{ij}(l(t))
\end{equation}
for any edge $\{ij\}\in E$.
As $\lim_{t\rightarrow +\infty}l(t)=l^*$, there exists $\xi_n\in (n, n+1)$ such that
$$\frac{dl_{ij}}{dt}|_{t=\xi_n}=l_{ij}(n+1)-l_{ij}(n)\rightarrow 0,$$
as $n\rightarrow +\infty$.
By the extended combinatorial Ricci flow (\ref{extended flow}), we have
$$\widetilde{K}_{ij}(l(\xi_n))=\frac{dl_{ij}}{dt}|_{t=\xi_n}\rightarrow 0.$$
Combining with (\ref{equ limit K_ij}), we have $\widetilde{K}_{ij}(l^*)=0$.
\qed

We have the following generalization of Theorem \ref{main theorem on convergence} for the extended combinatorial Ricci flow (\ref{extended flow}).

\begin{theorem}\label{context theorem convergence}
Suppose $(M, \mathcal{T})$ is an ideally triangulated compact $3$-manifold with boundary composed of surfaces of
negative Euler characteristic.
If  $(M, \mathcal{T})$ admits a hyper-ideal polyhedral metric $l^*\in \mathcal{L}(M, \mathcal{T})$ with zero combinatorial Ricci curvature,
then for any initial generalized hyper-ideal polyhedral metric $l(0)\in \mathbb{R}^E$ on $(M, \mathcal{T})$,
the solution $l(t)$ of the extended combinatorial Ricci flow (\ref{extended flow}) converges exponentially fast to $l^*$.
\end{theorem}
\proof
Recall that the function
$H$ in (\ref{energy function}) could be extended to be defined on $\mathbb{R}^E$ by
$$\widetilde{H}(l)=\sum_{\sigma\in T}\widetilde{F}_{\sigma}-2\pi\sum_{\{ij\}\in E}l_{ij},$$
which is a $C^1$ smooth convex function on $\mathbb{R}^E$ with $\nabla \widetilde{H}=-\widetilde{K}$ by Theorem \ref{extended co-volume}.
By the condition that $K(l^*)=0$ for $l^*\in \mathcal{L}(M, \mathcal{T})$, we have $\nabla \widetilde{H}(l^*)=0$ and
$l^*$ is a minimal point of $\widetilde{H}$ on $\mathbb{R}^E$.
Furthermore, $\lim_{l\rightarrow \infty}\widetilde{H}(l)=+\infty$ by the following property of convex functions, the proof of
which could be found in \cite{GX3} (Lemma 4.6).
This implies that $\widetilde{H}$ is a proper function on $\mathbb{R}^E$.

\begin{lemma}
Suppose $f(x)$ is a $C^1$ smooth convex function on $\mathbb{R}^n$ with $\nabla f(x_0)=0$ for some $x_0\in \mathbb{R}^n$,
$f(x)$ is $C^2$ smooth and strictly convex
in a neighborhood of $x_0$, then $\lim_{x\rightarrow \infty}f(x)=+\infty$.
\end{lemma}

Set $\phi(t)=\widetilde{H}(l(t))$, where $l(t)$ is a solution of the extended combinatorial Ricci flow (\ref{extended flow}).
Then $\phi(t)$ is a $C^1$ smooth function for $t\in [0, +\infty)$ with
\begin{equation}\label{monotonicity}
\frac{d\phi}{dt}=\nabla_l \widetilde{H}\cdot \frac{dl}{dt}=-\sum_{\{ij\}\in E}\widetilde{K}_{ij}^2\leq 0,
\end{equation}
which implies $\phi(t)$ is uniformly bounded for $t\in [0, +\infty)$, i.e.
$\widetilde{H}$ is bounded along the solution $l(t)$ of the extended combinatorial
Ricci flow (\ref{extended flow}).
By the properness of $\widetilde{H}$ on $\mathbb{R}^E$, we have the solution $l(t), t\in [0, +\infty)$, lies in a compact subset of $\mathbb{R}^E$.

By the boundedness of $\phi(t)$ on $[0, +\infty)$ and monotonicity of $\phi(t)$ from  (\ref{monotonicity}), we have
$\lim_{t\rightarrow +\infty} \phi(t)$ exists.
Therefore, there exists $\xi_n\in (n, n+1)$ such that
\begin{equation}\label{proof equ 1}
\phi(n+1)-\phi(n)=\frac{d\phi}{dt}|_{t=\xi_n}=-\sum_{\{ij\}\in E}\widetilde{K}_{ij}^2(l(\xi_n)) \rightarrow 0,
\end{equation}
as $n\rightarrow +\infty$.
By the boundedness of $l(t)$ for $t\in [0, +\infty)$, there exists a subsequence of $\xi_n$, still denoted by $\xi_n$ for simplicity,
such that $l(\xi_n)\rightarrow \overline{l}$ for some $\overline{l}\in \mathbb{R}^E$ as $n\rightarrow +\infty$,
which implies $\widetilde{K}(\overline{l})=0$ by (\ref{proof equ 1}) and the continuity of the generalized combinatorial Ricci curvature $\widetilde{K}$.

Note that $\widetilde{K}(\overline{l})=0$ for $\overline{l}\in \mathbb{R}^E$ and $K(l^*)=0$ for $l^*\in \mathcal{L}(M, \mathcal{T})$,
we have $\overline{l}=l^*$ by Theorem \ref{generalization of Luo-Yang rigidity}.
Therefore, there is a sequence $\xi_n\rightarrow +\infty$ such that $\lim_{n\rightarrow +\infty} l(\xi_n)=l^*$,
which implies that for $n$ large enough, $l(\xi_n)$ will lie in a small enough neighborhood of $l^*$.

Note that $l^*$ is a local attractor of the extended combinatorial Ricci flow (\ref{extended flow}) by Theorem \ref{negative definiteness of Lambda},
which was also observed by Luo \cite{L2},
we have the solution $l(t)$ of the flow (\ref{extended flow}) converges exponentially fast to $l^*$ by Lyapunov stability theorem (\cite{P}, Chapter 5).
\qed

\begin{remark}
Combinatorial Calabi flow for hyper-ideal polyhedral metrics
on ideally triangulated compact $3$-manifolds with boundary was introduced by Ge-Xu-Zhang \cite{GXZ},
which was defined as
\begin{equation}\label{combinatorial Calabi flow}
\frac{dl_{ij}}{dt}=-\Delta K_{ij},
\end{equation}
where $\Delta=\Lambda=(\frac{\partial K_{ij}}{\partial l_{kh}})$
is the combinatorial Laplace operator for hyper-ideal polyhedral metrics introduced by Luo \cite{L2}.
Combinatorial Calabi flow (\ref{combinatorial Calabi flow})
is a negative gradient flow of the combinatorial Calabi energy $\mathcal{C}(l)=\sum_{\{ij\}\in E} K_{ij}^2$ \cite{GXZ}.
As $\widetilde{a}_{ij}$ is $C^0$ and not $C^1$ by Remark \ref{remark non C1 of a_ij},
the combinatorial Laplace operator $\Delta=(\frac{\partial K_{ij}}{\partial l_{kh}})$
can not be extended continuously to be defined on $\mathbb{R}^E$
by Luo-Yang's extension \cite{LY}.
This causes that the extension method used for combinatorial Ricci flow (\ref{combinatorial Ricci flow}) in this paper
does not work for the combinatorial Calabi flow (\ref{combinatorial Calabi flow}).
Similar phenomenons happen for combinatorial Calabi flow for vertex scaling on closed surfaces \cite{ G, ZX}.
The convergence of the combinatorial Calabi flow for vertex scaling was proved in \cite{ZX}
using a different extension of combinatorial curvature introduced
in \cite{GGLSW,GLSW}. A key point of the proof in \cite{ZX} is that
the extension of combinatorial curvature introduced in \cite{GGLSW,GLSW} is $C^1$ smooth.
\end{remark}

As an application of Theorem \ref{context theorem convergence}, we have the following result,
which is a generalization of Corollary \ref{main corollary for regular triangulation}.

\begin{corollary}\label{corollary context}
Suppose $(M, \mathcal{T})$ is an ideally triangulated compact $3$-manifold with boundary.
\begin{description}
  \item[(1)] If $d_e\leq 6$ for any $e\in E$,
  there exists no hyper-ideal polyhedral metric in $\mathcal{L}(M, \mathcal{T})$ with zero combinatorial Ricci curvature.
  \item[(2)] If $d_e=N>6, \forall e\in E$, for some constant $N$, the solution of extended combinatorial Ricci flow (\ref{extended flow})
  converges to a hyper-ideal polyhedral metric in $\mathcal{L}(M, \mathcal{T})$ with zero combinatorial Ricci curvature
  for any initial generalized hyper-ideal polyhedral metric in $\mathbb{R}^E$.
\end{description}
\end{corollary}

\proof
For the first part of the corollary,
suppose there exists a hyper-ideal polyhedral metric $l^*\in \mathcal{L}(M, \mathcal{T})$
with zero combinatorial Ricci curvature.
Suppose $S$ is a component of the boundary $\partial M$, $|V|, |E|, |F|$ represent the number of vertices, edges and faces of $S$.
By the Gauss-Bonnet formula, we have $|V|-|E|+|F|=\chi(S)<0$.
Note that
$3|F|=2|E|$ and $2|E|=\sum_{x\in V}d_x\leq 6|V|$ by the condition $d_e\leq 6$, we have
$$0=\frac{1}{3}|E|-|E|+\frac{2}{3}|E|\leq |V|-|E|+|F|=\chi(S)<0,$$
which is impossible.
Therefore, there exists no hyper-ideal polyhedral metric in $\mathcal{L}(M, \mathcal{T})$ with zero combinatorial Ricci curvature in this case.

For the second part of the corollary, we just need to prove that there exists a hyper-ideal polyhedral metric $l^*\in \mathcal{L}(M, \mathcal{T})$
with zero combinatorial Ricci curvature, then the result follows from Theorem \ref{context theorem convergence}.
As the numbers of tetrahedra adjacent to the edges in $E$ are all the same,
we can assign each edge in $E$ with the same length $s\in (0, +\infty)$,
which corresponds to a hyper-ideal polyhedral metric in $\mathcal{L}(M, \mathcal{T})$.
In this case,  the lengths of vertex edges $x^i_{jk}$ are all the same with
$$\cosh x^i_{jk}=\frac{\cosh l_{ij}\cosh l_{ik}+\cosh l_{jk}}{\sinh l_{ij}\sinh l_{ik}}=\frac{\cosh s}{\cosh s-1},$$
which implies that
\begin{equation*}
\begin{aligned}
\cos a_{ij}
=\frac{\cosh x^i_{jk}\cosh x^i_{jh}-\cosh x^i_{kh}}{\sinh x^i_{jk}\sinh x^i_{jh}}=\frac{\cosh s}{2\cosh s-1}.
\end{aligned}
\end{equation*}
Therefore, the combinatorial Ricci curvature along any edge $\{ij\}\in E$ is
$$K(s)=2\pi-N\arccos \frac{\cosh s}{2\cosh s-1},$$
which is a strictly decreasing function of $s\in (0, +\infty)$.
Note that
$\lim_{s\rightarrow 0}K(s)=2\pi>0$
and
$\lim_{s\rightarrow +\infty}K(s)=\frac{\pi}{3}(6-N)<0$
by the condition $N> 6$. By the intermediate value theorem, there exists $s^*\in (0, +\infty)$ such that $K(s^*)=0$.
Therefore, $l^*=s^*(1,\cdots,1)^T$ is a hyper-ideal polyhedral metric in $\mathcal{L}(M, \mathcal{T})$ with zero combinatorial Ricci curvature.
\qed

\begin{remark}
By the proof of  Corollary \ref{corollary context} (1),
the condition on the number of tetrahedra adjacent to an edge in $E$ could be changed
to be the condition $d_x\leq 6$ on the degree of any vertex $x$ on the boundary surfaces.
Furthermore, for the nonexistence result in Corollary \ref{corollary context} (1),
we just need one component of the boundary satisfies the condition $d_x\leq 6$ for any vertex $x$
in the boundary component.
\end{remark}

\section{combinatorial Ricci flow for prescribed combinatorial Ricci curvature}\label{section 6}

We can modify the combinatorial Ricci flow to find hyper-ideal polyhedral metrics
with prescribed combinatorial Ricci curvature $\overline{K}$ as follows
\begin{equation}\label{modified flow}
\frac{dl_{ij}}{dt}=K_{ij}-\overline{K}_{ij}.
\end{equation}
Using the generalized combinatorial Ricci curvature $\widetilde{K}_{ij}$,
we can also extend the modified combinatorial Ricci flow (\ref{modified flow}) to the following form
\begin{equation}\label{extended modified flow}
\frac{dl_{ij}}{dt}=\widetilde{K}_{ij}-\overline{K}_{ij},
\end{equation}
which is called the extended modified combinatorial Ricci flow in the following.
The results for the modified combinatorial Ricci flow (\ref{modified flow}) and extended modified combinatorial Ricci flow (\ref{extended modified flow})
are paralleling to those of the combinatorial Ricci flow (\ref{combinatorial Ricci flow}) and the extended combinatorial Ricci flow (\ref{extended flow}).
We state the results as follows.
\begin{theorem}\label{theorem modified flow long time existence}
Suppose $(M, \mathcal{T})$ is an ideally triangulated compact $3$-manifold with boundary.
The solution of the extended modified combinatorial Ricci flow (\ref{extended modified flow}) exists for all time
for any initial generalized hyper-ideal polyhedral metric in $\mathbb{R}^E$
and uniquely extends the solution of the modified combinatorial Ricci flow (\ref{modified flow})
for any initial hyper-ideal polyhedral metric in $\mathcal{L}(M, \mathcal{T})$.
\end{theorem}

\begin{theorem}\label{theorem modified flow convergence}
Suppose $(M, \mathcal{T})$ is an ideally triangulated compact $3$-manifold with boundary.
If  $(M, \mathcal{T})$ admits a hyper-ideal polyhedral metric $\overline{l}\in \mathcal{L}(M, \mathcal{T})$
with combinatorial Ricci curvature $\overline{K}$,
then the solution $l(t)$ of the extended modified combinatorial
Ricci flow (\ref{extended modified flow}) converges exponentially fast to $\overline{l}$
for any initial generalized hyper-ideal polyhedral metric in $\mathbb{R}^E$.
\end{theorem}

The proof of Theorem \ref{theorem modified flow long time existence} is paralleling to that of Theorem \ref{main theorem on uniqueness and long time existence}
and the proof of Theorem \ref{theorem modified flow convergence} is similar to that of Theorem \ref{main theorem on convergence}
with the function $H(l)$ replaced by
$$\overline{H}(l)=\sum_{\sigma\in T}F_{\sigma}-2\pi\sum_{\{ij\}\in E}l_{ij}-\sum_{\{ij\}\in E}\overline{K}_{ij}l_{ij}$$
and the extension $\widetilde{H}(l)$ replaced by
$$\widetilde{\overline{H}}(l)=\sum_{\sigma\in T}\widetilde{F}_{\sigma}-2\pi\sum_{\{ij\}\in E}l_{ij}-\sum_{\{ij\}\in E}\overline{K}_{ij}l_{ij}.$$
As the proofs for Theorem \ref{theorem modified flow long time existence} and Theorem \ref{theorem modified flow convergence}
are almost the same as those of
Theorem \ref{main theorem on uniqueness and long time existence} and Theorem \ref{main theorem on convergence}
for the combinatorial Ricci flow (\ref{combinatorial Ricci flow})
and the extended combinatorial Ricci flow (\ref{extended flow}),
we omit the details of the proofs here.

\section{Some remarks and questions}\label{section 7}
In this paper, we have proved the global convergence of the extended combinatorial Ricci flow (\ref{extended flow})
to a complete hyperbolic metric on $(M, \mathcal{T})$ with totally geodesic boundary
under the condition that there exists a hyper-ideal polyhedral metric with zero combinatorial Ricci
curvature on the ideally triangulated compact 3-manifold with boundary $(M, \mathcal{T})$.
However, there may exist no hyper-ideal polyhedral metric with zero combinatorial Ricci curvature
on an ideally triangulated compact 3-manifold with boundary $(M, \mathcal{T})$,
while there exists a complete hyperbolic metric on the compact 3-manifold $M$ with totally geodesic boundary.
In other words, there may exist combinatorial obstacles for
the existence of hyper-ideal polyhedral metrics with zero combinatorial Ricci curvature on the ideally triangulated compact
3-manifold with boundary $(M, \mathcal{T})$.

An example of similar case is the vertex scaling for piecewise linear metric on surfaces \cite{L1}, for which there are combinatorial
obstacles for the existence of conformal factors with constant combinatorial curvature on a surface with a fixed triangulation.
For surfaces with fixed triangulations, Ge-Jiang \cite{GJ0} used the extension introduced by  Bobenko-Pinkall-Springborn \cite{BPS}
and Luo \cite{L3} to extend Luo's combinatorial Yamabe flow for vertex scaling on a triangulated surface \cite{L1} and proved the
convergence of the extended combinatorial Yamabe flow under the existence of a conformal factor with constant combinatorial
curvature on the surface with a fixed triangulation. The convergence of the extended combinatorial Yamabe flow
depends on the triangulation of the surface and the existence of
conformal factors with constant combinatorial curvature.

However, if one takes the piecewise constant curvature metric as a cone metric on the surface and
does surgery on the Delaunay triangulations of the surface by edge flipping,
it has been proved that there always exists a conformal factor defined on the vertices such that the
induced polyhedral metric on the surface has constant combinatorial curvature \cite{GGLSW,GLSW}.
Furthermore,
one can use the combinatorial Yamabe flow with surgery \cite{GGLSW,GLSW} and combinatorial Calabi flow with surgery \cite{ZX}
to find polyhedral metrics with constant combinatorial curvature, the convergence of which do not depend
on the initial triangulation of the surface and existence of conformal factors with constant combinatorial curvature.

For compact 3-manifolds with boundary, it is natural to ask the following question.

\textbf{Question:} Is there any way to do surgery on ideally triangulated compact 3-manifolds with boundary to ensure
the long time existence and convergence of the combinatorial Ricci flow (\ref{combinatorial Ricci flow})?

This is similar to the case of combinatorial Yamabe flow and combinatorial Calabi flow for vertex scaling on closed surfaces
and similar to the case of Hamilton's Ricci flow on 3-dimensional Riemannian manifolds.
Similar questions were also asked by Luo in \cite{L2}.
If one can give an affirmative answer to this question, it is conceive that one can give a solution to Luo's conjecture
without the assumption of existence of a hyper-ideal polyhedral metric with zero combinatorial Ricci curvature,
a proof of the convergence of combinatorial Calabi flow (\ref{combinatorial Calabi flow}) and
a new proof of Thurston's geometrization theorem for compact 3-manifolds
with boundary consisting of surfaces of negative Euler characteristic.
\\

\textbf{Acknowledgements}\\[8pt]
Part of the series of work, including \cite{X3,X4} and this paper, was done when the author was visiting the Rutgers University.
The author thanks Professor Feng Luo for his invitation to Rutgers University and communications and interesting on these work.
The author thanks Professor Feng Luo and Professor Tian Yang for explaining the details of their joint work \cite{LY} to the author.
Part of the series of work was reported in an online seminar of Rutgers University in April 2020.
The author thanks the participants in the seminar for communication and suggestions.
The research of the author is supported by
Fundamental Research Funds for the Central Universities and
National Natural Science Foundation of China under grant no. 61772379.

(Xu Xu) School of Mathematics and Statistics, Wuhan University, Wuhan 430072, P.R. China

E-mail: xuxu2@whu.edu.cn\\[2pt]


\begin{thebibliography}{50}
\setlength{\itemsep}{-1pt} \small

\bibitem{B} X. Bao, \emph{Hyperideal polyhedra in hyperbolic $3$-space}, Doctoral dissertation, University of Southern California, Los Angeles, 1998.

\bibitem{BB} X. Bao, F. Bonahon, \emph{Hyperideal polyhedra in hyperbolic $3$-space}, Bull. Soc. Math. France 130 (2002) 457-491.

\bibitem{BPS} A. Bobenko, U. Pinkall, B. Springborn, \emph{Discrete conformal maps and ideal hyperbolic polyhedra}.  Geom. Topol. 19 (2015), no. 4, 2155-2215.

\bibitem{B1} F. Bonahon, \emph{A Schl\"{a}fli-type formula for convex cores of hyperbolic $3$-manifolds}, J. Differential Geom. 50 (1998), 25-58.

\bibitem{CL} B. Chow, F. Luo, \emph{Combinatorial Ricci flows on surfaces}, J. Differential Geometry, 63 (2003), 97-129.


\bibitem{DG} S. Dai, H. Ge, \emph{Discrete Yamabe flows with R-curvature revisited}. J. Math. Anal. Appl. 484 (2020), no. 1, 123681, 11 pp.

\bibitem{FP} R. Frigerio, C. Petronio, \emph{Construction and recognition of hyperbolic $3$-manifolds with geodesic boundary}, Trans. AMS. 356 (2004), no. 8, 3243-3282.

\bibitem{F} M. Fujii, \emph{Hyperbolic $3$-manifolds with totally geodesic boundary}. Osaka J. Math. 27 (1990), no. 3, 539-553.

\bibitem{G} H. Ge, \emph{Combinatorial methods and geometric equations}, Thesis (Ph.D.)-Peking University, Beijing. 2012.

\bibitem{GH} H. Ge, B. Hua, \emph{$3$-dimensional combinatorial Yamabe flow in hyperbolic background geometry}, \href{https://arxiv.org/abs/1805.10643} {arXiv:1805.10643 [math.DG].}

\bibitem{GJ0} H. Ge, W. Jiang, \emph{On the deformation of discrete conformal factors on surfaces}. Calc. Var. Partial Differential Equations  55  (2016),  no. 6, Art. 136, 14 pp.

\bibitem{GJ1} H. Ge, W. Jiang, \emph{On the deformation of inversive distance circle packings, I}. Trans. Amer. Math. Soc. 372 (2019), no. 9, 6231-6261.

\bibitem{GJ2} H. Ge, W. Jiang, \emph{On the deformation of inversive distance circle packings, II}. J. Funct. Anal.  272  (2017),  no. 9, 3573-3595.

\bibitem{GJ3} H. Ge, W. Jiang, \emph{On the deformation of inversive distance circle packings, III}. J. Funct. Anal.  272  (2017),  no. 9, 3596-3609.

\bibitem{GJS} H. Ge, W. Jiang, L. Shen, \emph{On the deformation of ball packings}, \href{https://arxiv.org/abs/1805.10573} {arXiv:1805.10573 [math.DG].}

\bibitem{GM} H. Ge, S. Ma,  \emph{Discrete $\alpha$-Yamabe flow in $3$-dimension}. Front. Math. China 12 (2017), no. 4, 843-858.

\bibitem{GX1} H. Ge, X. Xu, \emph{Discrete quasi-Einstein metrics and combinatorial curvature flows in $3$-dimension},  Adv. Math. 267 (2014), 470-497.

\bibitem{GX2}  H. Ge, X. Xu, \emph{$\alpha$-curvatures and $\alpha$-flows on low dimensional triangulated manifolds}. Calc. Var. Partial Differential Equations 55 (2016), no. 1, Art. 12, 16 pp.

\bibitem{GX3}  H. Ge, X. Xu, \emph{On a combinatorial curvature for surfaces with inversive distance circle packing metrics}. J. Funct. Anal. 275 (2018), no. 3, 523-558.

\bibitem{GX4}  H. Ge, X. Xu, \emph{A combinatorial Yamabe problem on two and three dimensional manifolds}. \href{http://arxiv.org/abs/1504.05814}{arXiv:1504.05814v2 [math.DG].}

\bibitem{GXZ}  H. Ge, X. Xu, S. Zhang, \emph{Three-dimensional discrete curvature flows and discrete Einstein metrics}. Pacific J. Math. 287 (2017), no. 1, 49-70.

\bibitem{G1} D. Glickenstein, \emph{A combinatorial Yamabe flow in three dimensions}, Topology 44 (2005), No. 4, 791-808.

\bibitem{G2} D. Glickenstein, \emph{A maximum principle for combinatorial Yamabe flow}, Topology 44 (2005), No. 4, 809-825.

\bibitem{GGLSW} X. D. Gu, R. Guo, F. Luo, J. Sun, T. Wu, \emph{A discrete uniformization theorem for polyhedral surfaces II}, J. Differential Geom. 109 (2018), no. 3, 431-466.

\bibitem{GLSW} X. D. Gu, F. Luo, J. Sun, T. Wu, \emph{A discrete uniformization theorem for polyhedral surfaces}, J. Differential Geom.  109  (2018),  no. 2, 223-256.

\bibitem{H} R. S. Hamilton, \emph{Three-manifolds with positive Ricci curvature}. J. Differential Geom. 17 (1982), no. 2, 255-306.

\bibitem{Hart} P. Hartman, \emph{Ordinary differential equations}. John Wiley \& Sons, Inc., New York-London-Sydney 1964 xiv+612 pp.

\bibitem{L1} F. Luo, \emph{Combinatorial Yamabe flow on surfaces}, Commun. Contemp. Math. 6 (2004), no. 5, 765-780.

\bibitem{L2} F. Luo, \emph{A combinatorial curvature flow for compact $3$-manifolds with boundary}. Electron. Res. Announc. Amer. Math. Soc. 11 (2005), 12-20.

\bibitem{L2a} F. Luo, \emph{$3$-dimensional Schl\"{a}fli formula and its generalization}, Commun. Contemp. Math. 10 (2008), 835-842.

\bibitem{L3} F. Luo, \emph{Rigidity of polyhedral surfaces, III}, Geom. Topol. 15 (2011), 2299-2319.

\bibitem{LY} F. Luo, T. Yang,  \emph{Volume and rigidity of hyperbolic polyhedral 3-manifolds}. J. Topol. 11 (2018), no. 1, 1-29.

\bibitem{M} J. Milnor, \emph{The Schl\"{a}fli differential equality}, John Milnor Collected papers, Vol. 1, Geometry (Publish or Perish, Inc., Houston, TX, 1994).

\bibitem{Mo} E. Moise, \emph{Affine structures in $3$-manifolds V}. Ann. Math. (2) 56 (1952), 96-114.

\bibitem{P} L.S. Pontryagin, \emph{Ordinary Differential}, Addison-Wesley Publishing Company Inc., Reading, 1962.

\bibitem{S} J. Schlenker, \emph{Hyperideal polyhedra in hyperbolic manifolds}. \href{https://arxiv.org/abs/math/0212355v2}{arXiv:math/0212355v2.}

\bibitem{T1} W. Thurston, \emph{Geometry and topology of $3$-manifolds}, Princeton lecture notes 1976, \href{http://www.msri.org/publications/books/gt3m}{http://www.msri.org/publications/books/gt3m}.

\bibitem{V} E.B. Vinberg, \emph{Geometry. II}, Encyclopaedia of Mathematical Sciences, 29, Springer-Verlag, New York, 1988.

\bibitem{X1} X. Xu, \emph{Combinatorial $\alpha$-curvatures and $\alpha$-flows on polyhedral surfaces, I},  \href{https://arxiv.org/abs/1806.04516} {arXiv:1806.04516 [math.GT].}

\bibitem{X2} X. Xu, \emph{Combinatorial $\alpha$-curvatures and $\alpha$-flows on polyhedral surfaces, II}, preprint, 2018

\bibitem{X3} X. Xu, \emph{Combinatorial Ricci flow on cusped $3$-manifolds}, In preparation.

\bibitem{X4} X. Xu, \emph{Combinatorial Calabi flow on cusped $3$-manifolds}, In preparation.

\bibitem{Y} T. Yang, \emph{A combinatorial curvature flow for ideal triangulations}, Doctoral dissertation, University of Melbourne, 2019.

\bibitem{YJLG} X. Yin, M. Jin, F. Luo, David Gu, \emph{Computing constant-curvature metrics for hyperbolic $3$-manifolds with boundaries using truncated tetrahedral meshes}. Int. J. Shape Model. 14 (2008), no. 2, 169-188.

\bibitem{ZX} X. Zhu, X. Xu, \emph{Combinatorial Calabi flow with surgery on surfaces}. Calc. Var. Partial Differential Equations 58 (2019), no. 6, Paper No. 195, 20 pp.

\end{thebibliography}
\end{document}